 \documentclass[12pt]{article}
\usepackage[dvips]{epsfig}
\usepackage{amsthm,amsmath,amssymb}
\usepackage{latexsym}
\usepackage{placeins}
\usepackage{dblfloatfix}
\usepackage{afterpage}
\usepackage{color}
\usepackage{hyperref}
\hypersetup{
     colorlinks=true,
     linkcolor=blue,
     filecolor=blue,
     citecolor = black,      
     urlcolor=cyan,
     }
\definecolor{TWhite}{rgb}{1,1,1}
\definecolor{TBlack}{rgb}{0,0,0}
\definecolor{TBB}{rgb}{0,0,0.4}
\definecolor{LBlack}{rgb}{0.70,0.67,0.70}
\definecolor{Gray}{rgb}{0.55,0.52,0.55}
\definecolor{Comment}{rgb}{1,1,1}
\definecolor{TBlue}{rgb}{0,0,1}
\definecolor{LBlue}{rgb}{0.35,0.60,0.85}
\definecolor{TGreen}{rgb}{0,0.40,0.10}
\definecolor{LGreen}{rgb}{0.0,0.70,0.40}
\definecolor{Yellow}{rgb}{0.55,0.55,0}
\definecolor{TRed}{rgb}{1,0,0}
\definecolor{LRed}{rgb}{1.00,0.4,0.4}
\definecolor{TViolet}{rgb}{0.65,0,0.65}
\definecolor{LViolet}{rgb}{0.95,0.55,0.95}
\definecolor{TBrown}{rgb}{0.45,0.15,0}
\definecolor{LBrown}{rgb}{0.80,0.50,0}

\definecolor{blue}{rgb}{0,0,0.9}
\definecolor{red}{rgb}{0.9,0,0}
\definecolor{green}{rgb}{0,0.50,0.10}
\definecolor{violet}{rgb}{0.5804,0.0000,0.8275}

\newcommand{\violet}[1]{\begin{color}{violet}#1\end{color}}
\def\@themcountersep{}

\newtheorem{THEO}{Theorem}[section]
\newtheorem{ALGo}[THEO]{Algorithm}
\newtheorem{CONJ}[THEO]{Conjecture}
\newtheorem{COND}[THEO]{Condition}
\newtheorem{ASSUMP}[THEO]{Assumption}
\newtheorem{CORO}[THEO]{Corollary}
\newtheorem{DEFI}[THEO]{Definition}
\newtheorem{EXAMP}[THEO]{Example}
\newtheorem{FACT}[THEO]{Fact}
\newtheorem{HYPO}[THEO]{Hypothesis}
\newtheorem{LEMM}[THEO]{Lemma}
\newtheorem{PROB}[THEO]{Problem}
\newtheorem{PROP}[THEO]{Proposition}
\newtheorem{REMA}[THEO]{Remark}
\newcommand{\theo}{\begin{THEO}}
\newcommand{\algo}{\begin{ALGo} \rm}
\newcommand{\cond}{\begin{COND} \rm}
\newcommand{\assump}{\begin{ASSUMP} \rm}
\newcommand{\conj}{\begin{CONJ}}
\newcommand{\coro}{\begin{CORO}}
\newcommand{\defi}{\begin{DEFI} \rm}
\newcommand{\examp}{\begin{EXAMP} \rm}
\newcommand{\fact}{\begin{FACT}}
\newcommand{\hypo}{\begin{HYPO} \rm}
\newcommand{\lemm}{\begin{LEMM}}
\newcommand{\prob}{\begin{PROB} \rm}
\newcommand{\prop}{\begin{PROP}}
\newcommand{\rema}{\begin{REMA} \rm}
\newcommand{\etheo}{\end{THEO}}
\newcommand{\ealgo}{\end{ALGo}}
\newcommand{\econd}{\end{COND}}
\newcommand{\eassump}{\end{ASSUMP}}
\newcommand{\econj}{\end{CONJ}}
\newcommand{\ecoro}{\end{CORO}}
\newcommand{\edefi}{\end{DEFI}}
\newcommand{\eexamp}{\end{EXAMP}}
\newcommand{\efact}{\end{FACT}}
\newcommand{\ehypo}{\end{HYPO}}
\newcommand{\elemm}{\end{LEMM}}
\newcommand{\eprob}{\end{PROB}}
\newcommand{\eprop}{\end{PROP}}
\newcommand{\erema}{\end{REMA}}

\def\0{\mbox{\bf 0}}
\def\1{\mbox{\bf 1}}
\def\2{\mbox{\bf 2}}
\def\3{\mbox{\bf 3}}
\def\4{\mbox{\bf 4}}
\def\5{\mbox{\bf 5}}
\def\6{\mbox{\bf 6}}
\def\7{\mbox{\bf 7}}
\def\8{\mbox{\bf 8}}
\def\9{\mbox{\bf 9}}
\def\a{\mbox{\boldmath $a$}}
\def\b{\mbox{\boldmath $b$}}
\def\cc{\mbox{\boldmath $c$}}
\def\d{\mbox{\boldmath $d$}}

\def\h{\mbox{\boldmath $h$}}

\def\q{\mbox{\boldmath $q$}}

\def\u{\mbox{\boldmath $u$}}
\def\v{\mbox{\boldmath $v$}}
\def\w{\mbox{\boldmath $w$}}
\def\x{\mbox{\boldmath $x$}}
\def\y{\mbox{\boldmath $y$}}
\def\z{\mbox{\boldmath $z$}}

\def\AC{\mbox{$\cal A$}}

\def\inprod#1#2{\langle#1, \, #2\rangle}

\def\Real{\mbox{$\mathbb{R}$}}

\def\coneJ{\mbox{$\mathbb{J}$}}
\def\coneK{\mbox{$\mathbb{K}$}}
\def\coneL{\mbox{$\mathbb{L}$}}
\def\coneM{\mbox{$\mathbb{M}$}}
\def\coneN{\mbox{$\mathbb{N}$}}

\def\spaceE{\mbox{$\mathbb{E}$}}
\def\spaceV{\mbox{$\mathbb{V}$}}

\begin{document}

%\title{ \Large 
%Conic optimization problems with no duality gap
%} 

%\title{Strong duality of a conic optimization problem with two cones and a single equality 
%constraint}

\title{
Strong duality of a conic optimization problem with  a single hyperplane and two 
cone constraints
}

\author{
\normalsize 
Sunyoung Kim\thanks{Department of Mathematics, Ewha W. University, 52 Ewhayeodae-gil, Sudaemoon-gu, Seoul 03760, Korea 
			({\tt skim@ewha.ac.kr}). 
			 The research was supported  by   NRF 2021-R1A2C1003810.}, \and \normalsize
Masakazu Kojima\thanks{Department of Industrial and Systems Engineering,
	Chuo University, Tokyo 192-0393, Japan ({\tt kojima@is.titech.ac.jp}).
	This research was supported by Grant-in-Aid for Scientific Research (A) 19H00808 %26242027
%	and the Japan Science and 
%	Technology Agency (JST), the Core Research of Evolutionary Science and 
%	Technology (CREST) research project.             
 	}
%,
%  \and \normalsize
%Kim-Chuan Toh\thanks{Department of Mathematics, and Institute of Operations Research and Analytics, National University of Singapore,
%10 Lower Kent Ridge Road, Singapore 119076
%({\tt mattohkc@nus.edu.sg}). 
%This research is supported in part by the Ministry of Education, Singapore, Academic Research Fund (Grant number: R-146-000-257-112).
%} 
}

%\date{\normalsize Revised, \today}
\date{\normalsize \today}

\maketitle 
\vspace*{-0.4cm}

\begin{abstract}
\noindent
Strong (Lagrangian) duality of general conic optimization problems (COPs) has long been studied and its profound and  complicated results appear in different forms in a wide range of literatures.  As a result, characterizing the known and unknown results can  sometimes  be difficult. The aim of this article is to provide a unified and geometric view of strong duality of COPs for the known results. For our framework, we employ a COP minimizing a linear function in a vector variable $\x$ subject to a single hyperplane constraint $\x \in H$ and two cone constraints $\x \in \coneK_1$, $\x \in \coneK_2$.  It can be identically reformulated as a simpler COP  with the single hyperplane constraint $\x \in H$ and the single cone constraint $\x \in \coneK_1 \cap \coneK_2$. This simple COP and its dual as well as their duality relation can be represented geometrically, and  they have no duality gap without any constraint qualification.  The dual of the original target COP is equivalent to the dual of the reformulated COP if the Minkowski sum of the duals of the two cones $\coneK_1$ and $\coneK_2$ is closed or if the dual of the reformulated COP satisfies a certain Slater condition. Thus, these two conditions make it possible to transfer all duality results, including the existence and/or boundedness of optimal solutions,  on the reformulated COP to the ones on the original target COP, and further to the ones on a standard primal-dual pair of COPs with symmetry.

%Strong duality of general conic optimization problems (COPs) has long been studied and characterizing the known and unknown results can be difficult. We  provide a unified and geometric view of strong duality of COPs for the known results. We employ a COP minimizing a linear function in a vector variable $\x$ subject to a single hyperplane constraint $\x \in H$ and two cone constraints $\x \in \coneK_1$, $\x \in \coneK_2$.  It can be identically reformulated as a simpler COP with the single hyperplane constraint $\x \in H$ and the single cone constraint $\x \in \coneK_1 \cap \coneK_2$. This COP and its dual as well as their duality relation can be represented geometrically, and they have no duality gap without any constraint qualification. The dual of the original target COP is equivalent to the dual of the reformulated COP if the Minkowski sum of the duals of the two cones is closed or if the dual of the reformulated COP satisfies a Slater condition. With these two conditions, it is possible to transfer all duality results on the reformulated COP to the ones on the original target COP, and to the ones on a standard primal-dual pair of COPs with symmetry. 

\end{abstract}

\vspace{0.5cm}

\noindent
%{\bf Key words. } 
%Simple conic optimization problems, Newton method, Newton-bracketing method, quadratic convergence, conic relaxations, 
%Lagrangian DNN relaxations, linearly constrained QOPs in binary variables. 
{\bf Key words. } 
Duality,  conic optimization problems, simple conic optimization problems, 
generalizing  the Slater condition, 
closedness of the Minkowski sum of two cones

\bigskip

%\noindent 
%\red{
%\fbox{In my understanding, we can say "Slater's constraint" but neither "the Slater's constraint"}
%\fbox{nor "a Slater's constraint". So I have changed "Slater's" to "Slater" throughout the paper.}
%}
%Newton-bracketing method, quadratic convergence, 
%Lagrangian DNN relaxations, linearly constrained QOPs in binary variables. 

\vspace{0.5cm}

\noindent
{\bf AMS Classification.} 
90C20,  	%Quadratic programming 
90C22,  	%Semidefinite programming
90C25, 		%Convex programming
90C26,  	%Nonconvex programming, global optimization
90C27.		%Combinatorial optimization 

%!TEX root = ./main.tex

\section{Introduction}

It is well-known that  strong duality of conic optimization problems (COPs),  
including semidefinite programs (SDPs) \cite{ANJOS2012,NESTEROV1994,WOLKOWICZ2000},  
second order programs \cite{ALIZADEH2001,KIM2001}, doubly nonnegative programs 
\cite{ARIMA2018,YOSHISE2010} and 
copositive programs \cite{BURER11a,DUR2010}
depends on  their  representation.
In this paper, we consider a primal COP of the following form:
\begin{eqnarray}
\eta_p & = & 
\inf \; \left\{ \inprod{\q}{\x} : 
\x \in \coneK_1,  \ \x \in \coneK_2, \ \ \inprod{\h}{\x} = 1
\right\}.  \label{eq:COPprimal}
\end{eqnarray}
Here, 
\begin{eqnarray}
\left.
\begin{array}{lll} 
& & \spaceV  = \mbox{a finite dimensional vector space endowed with an inner product }  \\ 
& \ & \ \hspace{8mm} \mbox{$\inprod{\x}{\y}$ and a norm $\|\x\| = \sqrt{\inprod{\x}{\x}}$ for every $\x, \ \y \in \spaceV$},  \\ [2pt]
& & \coneK_1, \ \coneK_2 = \mbox{a nonempty closed convex cone in $\spaceV$}, \ 
\q \in \spaceV, \ \0\not=\h \in \coneK_1^*, \\[2pt]
&  & \coneK^* = \mbox{the dual } \{ \y \in \spaceV : \inprod{\x}{\y} \geq 0 \ \mbox{for every } \x \in \coneK\} \
\mbox{of a cone } \coneK \subset \spaceV. 
\end{array}
\right\} \label{eq:definition}
\end{eqnarray}
A general form of  COP with inequality constraints can be converted into COP~\eqref{eq:COPprimal}
 where the requirement 
$\0\not=\h \in \coneK_1^*$ above 
is naturally fulfilled. 
The conversion is described in  Section~\ref{section:symCOP}.

The (Lagrangian) dual of COP~\eqref{eq:COPprimal} is written as 
\begin{eqnarray}
\eta_d  & = & \sup \left\{ t_0 : \q - \h t_0 \in \coneK_1^* + \coneK_2^* \right\} 
 \label{eq:COPdual}\\ 
              & = & \sup \left\{ t_0 : \q - \h t_0 - \y_2 \in \coneK_1^* , \ \y_2 \in \coneK_2^* \right\},   \nonumber 
\end{eqnarray}
where 
$\coneK_1^* + \coneK_2^*$ denotes the Minkowski sum 
$\{ \y_1 + \y_2 : \y_1  \in \coneK_1^*,  \y_2 \in \coneK_2^* \}$  of $\coneK_1^*$ and $\coneK_2^*$. 
If primal COP~\eqref{eq:COPprimal} (or dual COP~\eqref{eq:COPdual}) is feasible, $\eta_p$ 
(or $\eta_d$) takes either a finite value or $-\infty$ (or $\infty$). If they are infeasible, 
we set $\eta_p = \infty$ and $\eta_d = -\infty$. Thus, the well-known weak duality 
$\eta_d \leq \eta_p$ holds even if one of them is infeasible. We simply write 
$-\infty < \eta_p < \infty$ (or $-\infty < \eta_d < \infty$) for the case where 
primal COP~\eqref{eq:COPprimal} (or dual COP ~\eqref{eq:COPdual}) is feasible and has 
a finite optimal value $\eta_p$ (or $\eta_d$). We say that {\em strong duality} 
holds if $-\infty < \eta_d = \eta_p < \infty$, and primal COP~\eqref{eq:COPprimal} (or dual COP ~\eqref{eq:COPdual}) is {\em solvable} if $\eta_p$ (or $\eta_d$) is attained by 
a feasible solution. 

Given a closed convex cone $\coneK \subset \spaceV$, 
primal COP~\eqref{eq:COPprimal} remains identical regardless of  how  
$\coneK$ is decomposed into the intersection $\coneK = \coneK_1\cap\coneK_2$ of two closed convex cones $\coneK_1$ and $\coneK_2$.
Dual COPs~\eqref{eq:COPdual} with different decompositions, however, may yield different duality results.
For instance, when we take $\coneK_1 = \coneK$ and $\coneK_2 = \spaceV$, 
there is no duality gap as we will see in Section~2. In the case where $\coneK_1 \not= \spaceV$ 
and $\coneK_2 \not= \spaceV$,  $\eta_d < \eta_p$ can occur. 

We introduce three conditions to characterize duality relations between 
COP~\eqref {eq:COPprimal} and COP~\eqref {eq:COPdual}
\begin{description}
\item{Cl:}
$\coneK_1^* + \coneK_2^*$ is closed. 
 \vspace{-2mm} 
\item{Ri:} 
There exists a $\widetilde{t} \in \Real$ such that 
$\widetilde{\y} \equiv \q- \h \widetilde{t}$ lies in 
$\mbox{relint}(\coneK_1^*+\coneK_2^*)$, 
the relative interior of 
$(\coneK_1^*+\coneK_2^*)$ with respect to the subspace spanned by 
$(\coneK_1^*+\coneK_2^*)$. 
\vspace{-2mm} 
\item{Po:} $\coneK_1\cap\coneK_2$ is pointed.
\vspace{-2mm}  
\end{description}

\noindent
It is known that the closure of $\coneK_1^* + \coneK_2^*$ coincides with $(\coneK_1\cap\coneK_2)^*$. 
Hence  Condition Cl  can be restated as 
$\coneK_1^* + \coneK_2^* = (\coneK_1\cap\coneK_2)^*$.  
% (see Lemma~\ref{lemma:MinkowskiSum}).
%Some Slater type sufficient conditions for 
%Conditions Cl, Ri and Po are given in Lemmas~\ref{lemma:threeCases} 
%and~\ref{lemma:sufficient}.  
We establish the following result. 

\theo \label{theorem:main0} Assume that $-\infty < \eta_p < \infty$ or $-\infty < \eta_d < \infty$ holds. 
\vspace{-2mm} 
\begin{description}
\item{(i)} If Condition {\rm Cl} or Condition {\rm Ri} is satisfied,  
then $-\infty < \eta_p = \eta_d < \infty$ holds.\vspace{-2mm}
\item{(ii)} If Condition {\rm Cl} 
is satisfied,  then dual {\rm COP}~\eqref{eq:COPdual} is solvable.
\vspace{-2mm}
\item{(iii)} The set of optimal  solutions of primal {\rm COP}~\eqref{eq:COPprimal} is nonempty and bounded 
if and only if Conditions {\rm Ri} and {\rm Po} are satisfied. 
\vspace{-2mm}
%\red{
\item{(iv)} 
The set of optimal  solutions of 
primal {\rm COP}~\eqref{eq:COPprimal} is nonempty and unbounded 
if {\rm Po} is not satisfied ({\it, i.e.}, 
$\coneK_1 \cap \coneK_2$ is not pointed) and Condition {\rm Ri} is satisfied.
%}
\vspace{-2mm} 
\end{description}
\etheo
\noindent 
Assertion (i) indicates that both Conditions Cl and Ri
are equally important 
to guarantee strong duality between primal COP~\eqref {eq:COPprimal} and dual COP~\eqref{eq:COPdual}. 
We note that ``Conditions Ri and Po'' in assertion 
(iii) can be replaced with 
a single condition that 
``there exists a $\widetilde{t} \in \Real$ such that 
$\q - \h \widetilde{t} \in \mbox{int}(\coneK_1^*+\coneK_2^*)$  
(the interior of $\coneK_1^*+\coneK_2^*$)''. 
(See Lemma~\ref{lemma:sufficient}).
The most significant (and straightforward) consequence of assertions (iii) 
and (i) of the theorem is:
\coro \label{coro:significant}
If the set of optimal solutions (or the feasible region) of primal  {\rm COP}~\eqref{eq:COPprimal} is 
nonempty and bounded, then $-\infty < \eta_p = \eta_d < \infty$.\vspace{-2mm}
\ecoro

\subsection*{Existing work and contribution of the paper}

Strong duality of COPs including SDPs has been widely studied  
\cite{AJAYI2020,LUO1997,NESTEROV1994,RAMANA1997,RAMANA1997a,ROCKAFELLAR1970,SHAPIRO2001}, 
and various results on their duality have been shown in different forms. 
In particular, Fenchel's duality theorem presented in Rockafellar \cite{ROCKAFELLAR1970} 
is quite general and (implicitly) covers a standard strong duality theorem of COPs under 
Slater's condition (\cite[Theorem 31.4]{ROCKAFELLAR1970}). 
Nesterov and Nemiroskii also gave a comprehensive discussion on duality for general 
COPs  in their book \cite{NESTEROV1994}, and presented a strong duality 
result under a Slater condition (see Theorem 4.2.1 of  \cite{NESTEROV1994}). 
Luo, Sturm and Zhang \cite{LUO1997} discussed \violet{the} boundedness of the optimal solution set 
of COPs in addition to their strong duality under Slater's condition. 
Shapiro proposed the closedness of a certain cone, instead of Slater condition, 
under which strong duality was established for a fairly general COPs 
\cite[Propositions 2.6 and 2.8]{SHAPIRO2001}. 
%As far as the authors are aware,  
%his closedness condition had been the weakest 
%one among the conditions to ensure strong duality. 
His closed condition has been playing a crucial role in the 
recent development of duality theory of COPs 
(see, \cite{AJAYI2020} and \cite{PATAKI2007}). 

The main purpose of this paper is to provide a  
unified and geometric overview of known strong duality results. 
A unique feature of the proposed approach in this paper is that strong duality of general COPs is obtained with a very simple form of COPs described geometrically by a cone, a 
hyperplane and a line. More precisely, we introduce the following two simple COPs: 
\begin{eqnarray*}
\mbox{P}(\coneK): & & \zeta_p(\coneK) = \inf\{ \inprod{\q}{\x} : \x \in \coneK, 
\ \inprod{\h}{\x} = 1 \}, \\ 
\mbox{D}(\coneJ): & & \zeta_d(\coneJ) = \sup\{ t : \q - \h t \in \coneJ \}. 
\end{eqnarray*}
Here $\coneK$ and $\coneJ$ are convex cones in a finite dimensional vector 
space $\spaceV$, $\q \in \spaceV$ and $\0 \not= \h \in \coneJ$. 
If we take $\coneJ = \coneK^*$, then COPs~\mbox{P}$(\coneK)$ and~\mbox{D}$(\coneK^*)$ 
serves as a primal-dual pair. 
The feasible region of primal COP P($\coneK$)  can be geometrically represented as 
the intersection of the hyperplane 
$H \equiv \{\x \in \spaceV : \inprod{\h}{\x} = 1 \}$ and the cone $\coneK$. 
Also the feasible region of dual COP D$(\coneK^*)$ can be viewed as the intersection of 
the $1$-dimension line $\{\q - \h t : t \in \Real \}$ with the dual cone $\coneK^*$ 
of $\coneK$, where $\h$ is a normal vector of the hyperplane $H$. 
This geometrical representation is an 
essential feature of the primal-dual pair of COPs P($\coneK$) 
and D$(\coneK^*)$, which makes it possible to geometrically interpret not only strong duality relation, but also duality gap
on the pair as we see in seven graphs in Figure 1. These  graphic interpretations are 
important contributions 
of this paper as  they show 
 basic essentials of strong duality geometrically.
In particular, strong duality 
$-\infty < \zeta_p(\coneK) = \zeta_d(\coneK^*) < \infty$ holds 
whenever either primal COP \mbox{P}$(\coneK)$ or dual COP \mbox{D}$(\coneK^*)$ has 
a finite optimal value {\em without any condition} (Theorem~\ref{theorem:simple}, (i)). 
If we choose $\coneK = \coneK_1 \cap \coneK_2$ and 
$\coneJ = \coneK_1^* + \coneK_2^*$, then the pair of 
COPs~\mbox{P}$(\coneK_1 \cap \coneK_2)$ 
and~\mbox{D}$(\coneK_1^* + \coneK_2^*)$  
represent the primal-dual pair of COPs~\eqref{eq:COPprimal} and~\eqref{eq:COPdual}.  
All assertions (i), (ii), (iii) and (iv) of Theorem~\ref{theorem:main0} on COPs~\eqref{eq:COPprimal} and~\eqref{eq:COPdual} follow  
from duality relations between 
COPs \mbox{P}$(\coneK)$ and~D$(\coneK^*)$.
Although they could be derived from some existing results referred above, 
the detailed derivation is not included here due to its  complexity.
Instead, we prefer to present 
a unified and geometric treatment of the strong duality through COPs P$(\coneK)$ 
and D$(\coneK^*)$. 

The simple primal-dual pair of COPs~\mbox{P}$(\coneK)$ and~\mbox{D}$(\coneK^*)$ 
was originally introduced in \cite{KIM2013} as a
Lagrangian-doubly nonnegative (DNN) relaxation of a class 
of quadratic optimization problems (QOPs), and their strong duality was 
studied in \cite[Lemma 2.5]{ARIMA2018}. 
Some relation of their results and Theorem~\ref{theorem:main0} will be discussed in 
Section~\ref{section:concludingRemarks}. 

\subsection*{Outline of the paper}

In Section~\ref{section:simpleCOP}, we discuss strong duality of the primal-dual pair of simple COPs P$(\coneK)$ and D$(\coneK^*)$, and establish Theorem~\ref{theorem:simple}. 
This theorem may be regarded as a special case of 
Theorem~\ref{theorem:main0} for $\coneK_1 = \coneK$ and $\coneK_2 = \spaceV$, 
and makes it easier to directly handle Theorem~\ref{theorem:main0}. 
Assertions of Theorem~\ref{theorem:simple} are illustrated in Figure 1.  
Based on Theorem~\ref{theorem:simple}, the main theorem of the paper, Theorem~\ref{theorem:main0}, is proved in Section~\ref{section:proofOfMainTheorem}. 
\violet{A} geometric implication of Conditions Cl and Ri is illustrated 
in Figure 2.
In addition, some lemmas  
are presented to relate Conditions Cl and Ri to Slater type conditions. 
In Section~\ref{section:symCOP}, we apply Theorem~\ref{theorem:main0} to 
a standard primal-dual pair of COPs with symmetry, and present a slight extension of 
the strong duality results given in \cite[Propositions 2.6 and 2.8]{SHAPIRO2001}. 
%Some remarks on \cite[Theorem 3.1]{AJAYI2020} related to Theorem~\ref{theorem:main0} 
%are also provided. 
We conclude in Section~5 with remarks on the strong duality result 
in \cite{KIM2013,ARIMA2018} for the pair of COPs of P$(\coneK)$ and 
D$(\coneK^*)$ which was
originally proposed as a Lagrangian-DNN relaxation of a class of QOPs. 
%Section~\ref{section:appendix} 
%Appendix includes supplementary materials on 
%technical details of some lemmas.

\subsection*{Notation and symbols}

For every convex subset $C$ of a finite dimensional space $\spaceV$, 
int$C$, linspan$C$, 
relint$C$, cl$C$ and $C^{\perp}$
denote the interior of $C$, the linear subspace spanned by $C$ 
({\it i.e.}, the smallest linear subspace containing $C$), 
the relative interior of $C$ with 
respect to linspan$C$, the closure of $C$ and the orthogonal complement of 
linspan$C$, respectively.

%!TEX root = ./main.tex

\section{Strong duality between COPs~\mbox{P}$(\coneK)$ and~\mbox{D}$(\coneK^*)$}

\label{section:simpleCOP}

In this section, we establish Theorem~\ref{theorem:simple} 
which leads to assertions (i), (ii), (iii) and (iv) of  our main theorem, Theorem~\ref{theorem:main0}. 
Throughout this section we assume that $\coneK$ is a closed convex cone in 
a finite dimensional vector space $\spaceV$, 
$\q \in \spaceV$ and $\0 \not= \h \in \coneK^*$.
To facilitate  connecting assertions of Theorem~\ref{theorem:simple} 
with those of 
Theorem~\ref{theorem:main0}, 
we modify 
Conditions Cl, Ri and Po with an arbitrary convex cone $\coneJ \subset \spaceV$. 
\vspace{-2mm}
\begin{description}
\item[{\normalfont Cl'$(\coneJ)$}:] $\coneJ$ is closed.\vspace{-2mm} 
\item[{\normalfont Ri'$(\coneJ)$}:] There exists a $\widetilde{t} \in \Real$ such that 
$\widetilde{\y} \equiv \q - \h \widetilde{t} \in {\rm relint}\coneJ$.
\vspace{-2mm} 
\item[{\normalfont Po'$(\coneJ)$}:] $\coneJ$ is pointed.\vspace{-2mm}
\end{description}

\noindent
Note that original Conditions Cl, Ri and Po can be rewritten as 
Cl'$(\coneK_1^*+\coneK_2^*)$, Ri'$(\coneK_1^*+\coneK_2^*)$ and 
Po'$(\coneK_1\cap\coneK_2)$, respectively. 

\theo \label{theorem:simple} 
Assume that 
% \begin{eqnarray}
$
-\infty < \zeta_p(\coneK) < \infty$ or  $-\infty < \zeta_d(\coneK^*) < \infty 
\label{eq:basicAssumption}
$ 
% \end{eqnarray}
holds (as in Figure 1 (a), (b), (c) or (e)).
\vspace{-2mm}
\begin{description}
\item{(i)} 
$-\infty < \zeta_p(\coneK) = \zeta_d(\coneK^*)  < \infty$ holds. 
\vspace{-2mm}
\item{(ii)} 
Dual {\rm COP D}$(\coneK^*)$  is solvable. \vspace{-2mm}
\item{(iii)} 
The set of optimal  solutions of primal 
{\rm COP} {\rm P}$(\coneK)$ is nonempty and bounded if and only if 
Conditions {\rm Ri'}$(\coneK^*)$ and {\rm Po'}$(\coneK)$
are satisfied (as in Figure 1 (a) and (e)). 
\vspace{-2mm}
%\red{
\item{(iv)} 
%In addition to \eqref{eq:basicAssumption},  
%assume that {\rm Po'(\coneK)} is not satisfied, {\it i.e.}, $\coneK$ is not 
%pointed. Then 
%
The set of optimal  solutions of primal 
{\rm COP} {\rm P}$(\coneK)$ is nonempty and unbounded if 
{\rm Po'(\coneK)} is not satisfied ({\it i.e.}, $\coneK$ is not 
pointed) and Condition {\rm Ri'}$(\coneK^*)$ is satisfied 
 (as in Figure 1 (c)). 
\vspace{-2mm} 
\end{description}
\etheo

The geometric interpretation of assertions (i), (ii), (iii) and (iv) 
is illustrated with seven two-dimensional examples in Figure 1 and Table 1 to help 
the reader understand them.
A proof of Theorem~\ref{theorem:simple}  is presented after
a basic lemma below.

\afterpage{
\begin{figure}[t!]  \vspace{-0.3cm} 
\begin{center}
\includegraphics[width=0.32\textwidth]{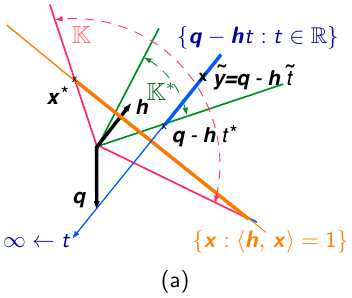}
\includegraphics[width=0.32\textwidth]{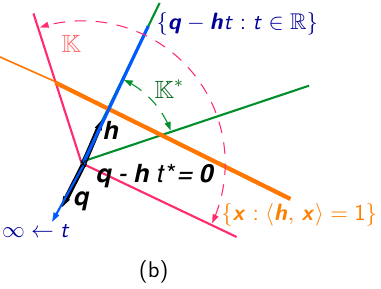}

\includegraphics[width=0.32\textwidth]{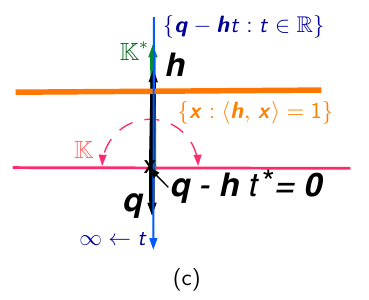}
\includegraphics[width=0.32\textwidth]{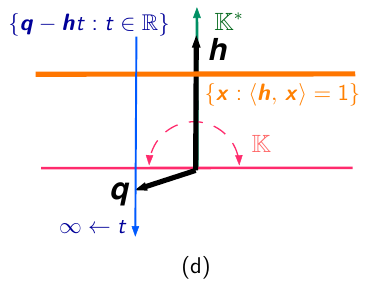}

\includegraphics[width=0.32\textwidth]{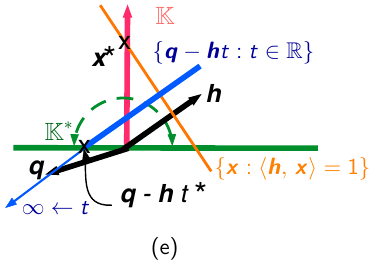}
\includegraphics[width=0.32\textwidth]{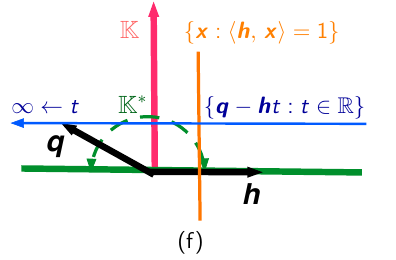}
\includegraphics[width=0.32\textwidth]{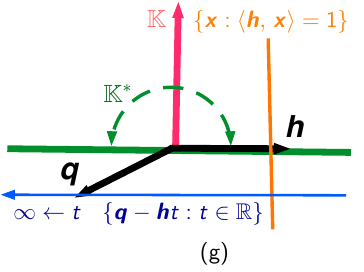}
\end{center}

\caption{ Cases (a), (b), (c) and (e) satisfy the 
assumption of Theorem~\ref{theorem:simple} (hence $\zeta_p(\coneK) = \zeta_d(\coneK^*)$ 
holds by Theorem~\ref{theorem:simple} (i) and D$(\coneK^*)$ is solvable), 
while cases (d), (f) and (g) do not. 
Cases (a) and (e) satisfy Conditions Ri'($\coneK^*$) and 
Po'($\coneK$); hence the set of optimal solutions 
of P($\coneK$)  is nonempty and bounded by Theorem~\ref{theorem:simple} (ii).
In case (b), Po'($\coneK$) is satisfied and 
both P($\coneK$) and D($\coneK^*$) are solvable, 
but Condition Ri'($\coneK^*$) is violated. Hence we know 
 that Condition Ri'($\coneK^*$) is merely a sufficient 
condition for the solvability of P($\coneK$) but not necessary. 
Case (b) is sensitive to the change in $\q$: 
if the direction of $\q$ is perturbed slightly, then  
either Ri'($\coneK^*$) together with Po'$(\coneK)$ is satisfied or D($\coneK^*$) gets infeasible (or, either P($\coneK$) has a unique optimal solution at $\x^*$ or $\zeta_p(\coneK) 
= -\infty$). Case (c) satisfies Condition Ri'($\coneK^*$) but does not 
satisfy Condition Po'$(\coneK)$; hence the set of solutions of 
P($\coneK$) is nonempty and unbounded.
Case (d) is an example for 
$- \infty  = \zeta_d(\coneK^*) = \zeta_p(\coneK)$, and case (f) for 
$\zeta_d(\coneK^*) = \zeta_p(\coneK) = \infty $. In case (g), both P($\coneK$)
and D($\coneK^*$) are infeasible; hence $-\infty = \zeta_d(\coneK^*) < \zeta_p(\coneK) = \infty$. Table~1 summarizes all cases, where cases designated by $\times$ cannot occur 
by weak duality or (i) of Theorem 2.1.
}
\end{figure} 
\begin{table}[h!]
	\begin{center}
%\scriptsize{
		\begin{tabular}{|cc||c|c|c|c|c|c|}
			\hline
                                            &                                   & \multicolumn{3}{c|}{$\mbox{P}(\coneK)$}      \\
                                            &                                   & \multicolumn{2}{c|}{Feasible} & Infeasible \\
                                            &                                   &  $-\infty <\zeta_p(\coneK)$ 
                                            & $-\infty =\zeta_p(\coneK)$ & $\zeta_p(\coneK) = \infty$\\
\hline \hline             
                                            & Feasible     $\zeta_d(\coneK^*) < \infty$ 
                                            & (a),(b),(c),(e) & $\times$(by Th.2.1 (i)) & $\times$(by Th.2.1 (i))\\
D$(\coneK^*)$ & Feasible     $ \zeta_d(\coneK^*) = \infty$ &$\times$(by Th.2.1 (i)) &$\times$(by w.duality) & (f)\\
                                            & Infeasible  $ \zeta_d(\coneK^*) = - \infty$ & $\times$(by Th.2.1 (i)) & (d) & (g) \\
\hline                                             
\end{tabular}	
%}
	\end{center}
	\caption{Cases (a) through (g) correspond to those of Figure 1, respectively. Cases designated by $\times$ cannot occur by weak duality or 
(i) of Theorem~\ref{theorem:simple}.}
\end{table}
}

\lemm \label{lemma:RiAndPo}
Let $\coneJ \subset \spaceV$ be a nonempty closed convex cone. Then Conditions 
{\rm Ri'}$(\coneJ^*)$ and {\rm Po'}$(\coneJ)$ hold if and only if 
\begin{eqnarray}
\mbox{there exists a $\widetilde{t} \in \Real$ such that } \ 
\widetilde{\y} \equiv \q - \h\widetilde{t} \in \mbox{int}(\coneJ^*). 
\label{eq:RiAndPo} 
\end{eqnarray}
\elemm 
\proof{We know that $\coneL \subset \coneJ$ holds for a linear subspace 
$\coneL$ if and only if $\coneJ^* \subset \coneL^L$. 
Hence $\coneJ$ is pointed if and only if int$(\coneJ^*) \not= \emptyset$. 
This implies the desired result. 
\qed
}

\smallskip

\noindent
{\em Proof of Theorem~\ref{theorem:simple} (i) and (ii):} Let $\d \in \mbox{relint}\coneK$. % Then there exists an $\alpha > 1$ such that 
It follows from $\h \in \coneK^*$ that $\inprod{\h}{\d} \geq 0$. Assume that 
$\inprod{\h}{\d} = 0$. Since $\d \in \mbox{relint} \coneK$, for every $\x \in \coneK$, 
there is a $\mu > 1$ such that $(1-\mu)\x +\mu\d \in \coneK$ 
(see \cite[Theorem 6.4]{ROCKAFELLAR1970}). Hence $\0 \leq \inprod{\h}{(1-\mu)\x +\mu\d}
=(1-\mu)\inprod{\h}{\x}$ for every $\x \in \coneK$, which together with $\h \in \coneK^*$ implies that $\h \in \coneK^{\perp}$. Therefore P$(\coneK)$ is infeasible or 
$\zeta_p(\coneK) = +\infty$, and either $\zeta_d(\coneK^*) = +\infty$ 
(if $\q \in \coneK^*$) or $\zeta_d(\coneK^*) = -\infty$ (if $\q \not\in \coneK^*$). 
This contradicts the assumption, and we have shown that $\inprod{\h}{\d} > 0$. 
Now, let $\tilde{\x} = \d / \inprod{\h}{\d}$. Then $\tilde{\x}$ is a feasible solution of 
P$(\coneK)$ satisfying Slater's condition, {\it i.e.} $\tilde{\x} \in \mbox{relint}\coneK$
 and $\inprod{\h}{\tilde{\x}} = 1$. By the standard strong duality theorem (for example, 
 \cite[Theorem 7]{LUO1997}, \cite[Theorem 31.4]{ROCKAFELLAR1970}), assertions (i) and (ii) hold. 
\qed

\smallskip

\noindent
{\em Proof of Theorem~\ref{theorem:simple} (iii):} 
We already by (ii) know that $-\infty < \zeta_p(\coneK)=\zeta_d)\coneK^*) < \infty$. Also, 
by Lemma~\ref{lemma:RiAndPo}, we may replace Conditions {\rm Ri'}$(\coneK^*)$ 
and {\rm Po'}$(\coneK)$ with~\eqref{eq:RiAndPo}. 

``if part'': 
We assume % on the contrary 
that \eqref{eq:RiAndPo} holds and  
on the contrary that the set of optimal solutions of 
P$(\coneK)$ is either empty or unbounded, and we show a contradiction. 
Since $-\infty < \zeta_p(\coneK) < \infty$ by  assumption, 
there exist an $\epsilon \geq 0$ and a sequence of feasible solutions $\{\x^k\}$ of 
P$(\coneK)$ such that 
\begin{eqnarray}
\x^k \in \coneK, \ \inprod{\q}{\x^k}\leq \zeta_p(\coneK) + \epsilon, 
\ \inprod{\h}{\x^k} = 1 \label{eq:sequence}
\end{eqnarray} 
for $k=1,2,\ldots$, and $\|\x^k\| \rightarrow \infty$  as $k \rightarrow \infty$. 
(If the set of optimal solutions is 
nonempty, we can take $\epsilon = 0$ and 
$\inprod{\q}{\x^k} = \zeta_p(\coneK)$).  
We may assume without 
loss of generality that $\x^k / \|\x^k\| \in \coneK$ converges to some 
$\Delta\x \in \coneK$ with $\|\Delta\x \| = 1$. 
Dividing the relation \eqref{eq:sequence}  by  $\|\x^k\|$ and taking the 
limit as $k \rightarrow \infty$, we see that
\begin{eqnarray*}
 \Delta\x \in \coneK, \ \|\Delta\x \| = 1, \ \inprod{\q}{\Delta\x} \leq 0, \ 
 \inprod{\h}{\Delta\x} = 0. 
\end{eqnarray*}
By \eqref{eq:RiAndPo}, 
if we take a sufficiently small $\delta > 0$ such that 
$ 
\widetilde{\y} -\delta\Delta\x = \q - \h \widetilde{t} -\delta\Delta\x \in \coneK^*,  
$
then 
$ 0 \leq \inprod{\q - \h\widetilde{t}- \delta\Delta\x}{\Delta\x}  
\leq - \delta < 0,$ 
which is a contradiction.

``only if part'': 
Assume on the contrary that 
the line  $S = \left\{\z = \q - \h t : t \in\Real \right\}$ does not intersect 
with int$(\coneK^*)$. 
By the separation theorem on convex sets (see, for example, \cite[Theorem 11.3]{ROCKAFELLAR1970}),  
there exist an $\alpha \in \Real$  and a nonzero $\Delta\x \in \spaceV$ such that 
\begin{eqnarray*}
\inprod{\y}{\Delta\x} \geq \alpha \geq \inprod{\z}{\Delta\x} \ \mbox{if } \y \in 
\coneK^* \ \mbox{and } \z \in S.
\end{eqnarray*}
Since $\coneK^*$ is a closed cone containing $\0 \in \spaceV$, we see that 
\begin{eqnarray*}
& & \inprod{\y}{\Delta\x} \geq 0 \ \mbox{for every } \y \in \coneK^*; \ 
\mbox{hence } \0 \not= \Delta\x \in \coneK, \\
& & 0 \geq \alpha \geq \inprod{\q - \h t}{\Delta\x} \ \mbox{for every } t \in \Real; \ \mbox{hence } \
 \inprod{\h}{\Delta\x} = 0 \ \mbox{and }  0  \geq \inprod{\q}{\Delta\x}.
\end{eqnarray*}
Let $\x^*$ be an optimal solution of P$(\coneK)$  whose existence is guaranteed by the assumption of the ``only if" part. Then, 
\begin{eqnarray*}
& & \x^* + \mu \Delta\x \in \coneK, \ \inprod{\h}{\x^*+\mu\Delta\x} = 1; \\
& & \mbox{(hence $\x^* + \mu \Delta\x$ is a feasible solution of P$(\coneK)$)}, \\
& & \zeta_p(\coneK) = \inprod{\q}{\x^*} \geq 
\inprod{\q}{\x^* + \mu\Delta\x} 
\geq \inprod{\q}{\x^*} = \zeta_p(\coneK)
\end{eqnarray*}
hold for all $\mu \geq 0$. This implies that $\{\x^*+\mu\Delta\x: 
\mu \geq 0 \}$ forms an unbounded ray in the set of optimal solutions of 
${\rm P}(\coneK)$.  
This contradicts to the assumption of ``only if" part. 
\qed

\smallskip

\noindent
{\em Proof of Theorem~\ref{theorem:simple} (iv):} 
Assume that $\coneK$ is not pointed. We know by (iii) that 
the set of optimal solution of P$(\coneK)$ is either empty or unbounded. 
Thus it suffices to show that the set of optimal solution is nonempty 
if Condition Ri'$(\coneK^*)$ is satisfied.  
To show this, we embed the primal COP P$(\coneK)$ 
to an equivalent COP P$(\coneJ)$ in a subspace $\coneL$ of $\spaceV$, 
\violet{to which} we will apply assertion (iii). We take the 
subspace of $\spaceV$ spanned by $\coneK^*$ for the subspace $\coneL$. 
Let $\Pi$ denote the orthogonal projection from $\spaceV$ 
onto $\coneL$, and define the convex cone $\coneJ = \{ \Pi(\x) : \x \in \coneK\}$. 
Then we observe that 
\begin{eqnarray*}
& & \coneJ \subset \coneL, \ \coneK^* \subset \coneL, \ \h \in \coneK^* \subset \coneL, \ \coneL^{\perp} \subset \coneK \ \mbox{and}  \\ 
& & \q \in \coneL  \ \mbox{(since otherwise $\zeta_p(\coneK) = -\infty$)}. 
\end{eqnarray*}
Hence we may regard $\coneL$ as the entire space on which the primal-dual pair 
P$(\coneJ)$ and D$(\coneJ^*)$ are defined, where $\coneJ^* = \{\y \in \coneL : 
\inprod{\x}{\y} \geq 0 \ \mbox{for every } \x \in \coneJ\} = \coneK^*$. 
Furthermore, for every feasible solution $\x$ of P$(\coneK)$, 
%  and $\delta\x \in \coneL^{\perp}$,  
\begin{eqnarray*}
& & 1 = \inprod{\h}{\x} = \inprod{\h}{\Pi(\x)}, \ \Pi(\x) \in \coneJ \
% & & \inprod{\q}{\delta\x} = 0 \ \mbox{(otherwise $\zeta_p(\coneK) = - \infty$}), \
\mbox{and } \inprod{\q}{\x} =  \inprod{\q}{\Pi(\x)}. 
\end{eqnarray*}
This implies the equivalence of P$(\coneK)$ to P$(\coneJ)$ in the sense that 
$\x$ is a feasible solution of P$(\coneK)$ if and only 
if its projection $\Pi(\x)$ onto $\coneL$ is a feasible solution of 
P$(\coneJ)$ with the same objective value 
$\inprod{\q}{\Pi(\x)} = \inprod{\q}{\x}$. 
By construction, $\coneJ$ is pointed with respect the space $\coneL$. 
By assertion (iii), the set of optimal solutions of P$(\coneJ)$ is nonempty and 
bounded if and only if Condition {\rm Ri'}$(\coneJ^*)$ 
(= Condition {\rm Ri'}$(\coneK^*)$) holds. Therefore, the set of optimal solutions 
of P$(\coneK)$ is nonempty if Condition {\rm Ri'}$(\coneK^*)$ holds. 
\qed

\rema \label{remark:(iii)} 
The dual version of (iii) was given in \cite[Theorem 5]{LUO1997}. 
\erema

%!TEX root = ./main.tex

\section{Proofs of Theorem~\ref{theorem:main0} and related results}

\label{section:proofOfMainTheorem}

Throughout this section, we use the notation and symbols given in~\eqref{eq:definition}.   
We need the following lemmas to prove Theorem~\ref{theorem:main0}. 

\lemm \label{lemma:MinkowskiSum} 
Let $\coneK_1$ and $\coneK_2$ be nonempty closed 
convex cones in $\spaceV$. Then $(\coneK_1\cap\coneK_2)^* = $ 
cl$(\coneK_1^* + \coneK_2^*)$ and 
% \red{
{\rm relint}$((\coneK_1\cap\coneK_2)^*) = $ {\rm relint}$(\coneK_1^* + \coneK_2^*)$. 
% }.
 \vspace{-2mm}
\elemm
\proof{
Assertion $(\coneK_1\cap\coneK_2)^* = $ cl$(\coneK_1^* + \coneK_2^*)$ is 
well-known (see, for example, \cite{PATAKI2007}). 
Assertion  relint$((\coneK_1\cap\coneK_2)^*) = $ relint$(\coneK_1^* + \coneK_2^*)$ follows from cl$(\coneK_1^* + \coneK_2^*) = 
 (\coneK_1 \cap \coneK_2)^*$ and the convexity of $\coneK_1^* + \coneK_2^*$.
\qed

\smallskip

\lemm \label{lemma:CoditionRi}
Let $\coneK = \coneK_1 \cap \coneK_2$. 
Assume that  $-\infty < \eta_d < \infty$ and 
Condition {\rm Ri'}$(\coneK_1^*+\coneK_2^*)$  is satisfied.
Then 
$\eta_d = \zeta_d(\coneK_1^*+\coneK_2^*) = \zeta_d(\coneK^*)$.\vspace{-2mm}
\elemm
\proof{
The identity $\eta_d = \zeta_d(\coneK_1^*+\coneK_2^*)$ holds by definition. 
Since $\zeta_d(\coneK_1^*+\coneK_2^*) \leq \zeta_d(\coneK^*)$ holds by  
$ \coneK_1^* +\coneK_2^* \subset \coneK^*$, it suffices to show that 
$\zeta_d(\coneK^*) \leq \zeta_d(\coneK_1^*+\coneK_2^*)$. 
Let $t$ be a feasible solution of dual COP D$(\coneK^*)$;  
$\q - \h t \in \coneK^* = $ cl$(\coneK_1^* + \coneK_2^*)$ by Lemma~\ref{lemma:MinkowskiSum}.  
By Condition {\rm Ri'}$(\coneK_1^*+\coneK_2^*)$, we know that 
$\q - \h \widetilde{t} \in $ relint$(\coneK_1^* + \coneK_2^*)$. 
Hence the convex combination $(1-\epsilon) (\q - \h  t) + \epsilon ( \q - \h \widetilde{t})$ 
with any small $\epsilon \in (0,1]$ lies in relint$(\coneK_1^* + \coneK_2^*)$. It follows that 
\begin{eqnarray*}
    \coneK_1^* + \coneK_2^* \supset  \mbox{relint}(\coneK_1^* + \coneK_2^*)  \ni  (1-\epsilon) (\q - \h  t) + \epsilon ( (\q - \h \widetilde{t}) = \q - \h (t  + \epsilon(\widetilde{t}-t))
\end{eqnarray*}
holds for any small $\epsilon \in (0,1]$.
This implies that $(t + \epsilon(\widetilde{t}-t))$ is a feasible solution of 
dual COP D$(\coneK_1^* +\coneK_2^*)$ for any small $\epsilon \in (0,1]$. 
Therefore we conclude that $\zeta_d(\coneK^*) \leq \zeta_d(\coneK_1^* +\coneK_2^*)$. 
\qed
} 

\smallskip

\noindent
{\it Proof of Theorem~\ref{theorem:main0} (i):} 
Let $\coneK = \coneK_1 \cap \coneK_2$. 
Then primal COP~\eqref{eq:COPprimal} and P$(\coneK)$ are identical  
(hence $\eta_p = \zeta_p(\coneK))$. If Condition Cl is satisfied, then 
dual COP~\eqref{eq:COPdual} and D$(\coneK^*)$ are identical 
(hence $\eta_d = \zeta_d(\coneK^*)$) by Lemma~\ref{lemma:MinkowskiSum},  
and if Condition Ri is satisfied, 
then $\eta_d = \zeta_d(\coneK^*)$ holds by Lemma~\ref{lemma:CoditionRi}. 
Therefore, assertion (i) follows from assertion (i) of Theorem~\ref{theorem:simple}. 

\smallskip

\noindent
\noindent
{\it Proof of Theorem~\ref{theorem:main0} (ii):} 
Let $\coneK = \coneK_1\cap  \coneK_2$. Then, by Lemma~\ref{lemma:CoditionRi}, 
primal COP~\eqref{eq:COPprimal} and 
dual COP~\eqref{eq:COPdual} are identical to P($\coneK$) and D$(\coneK^*)$, 
respectively. Hence 
assertion (ii) follows from  assertion (ii)  of  Theorem~\ref{theorem:simple}. 
\qed

\smallskip

\noindent
{\it Proof of Theorem~\ref{theorem:main0} (iii) and (iv):} 
Let $\coneK = \coneK_1 \cap \coneK_2$. 
Then we can replace primal COP~\eqref{eq:COPprimal} with P$(\coneK)$ 
in assertions (iii) and (iv) since they are identical. Also Po and Po'$(\coneK)$ 
are identical. 
By Lemma~\ref{lemma:MinkowskiSum}, we also see 
relint$(\coneK_1^*+\coneK_2^*) = $ relint$(\coneK^*)$.  
Hence Condition Ri'$(\coneK_1^*+\coneK_2^*)$ ($=$ Condition Ri) holds if 
and only if Condition Ri'$(\coneK^*)$ holds. 

``if part of (iii)'': 
Assume that Conditions Ri'$(\coneK^*)$ and Po'$(\coneK)$ are satisfied. Then 
$\zeta_d(\coneK_1^*+\coneK_2^*) = \zeta_d(\coneK^*)$ by Lemma~\ref{lemma:CoditionRi}. Thus, $-\infty < \zeta_p(\coneK) < 
\infty$ or $-\infty < \zeta_d(\coneK^*) < \infty$ by the assumption of 
the theorem. Therefore, the set of optimal solutions of P($\coneK$), which 
coincides with the set of optimal solutions of primal COP 
P($\coneK_1\cap\coneK_2$), is nonempty and bounded by (iii) of Theorem~\ref{theorem:simple}. 

``only if part of (iii)'': 
Conversely, assume that the set of optimal solutions of 
P($\coneK_1\cap\coneK_2$) $=$ P($\coneK$) is nonempty and bounded, which implies  
$-\infty < \zeta_p(\coneK) < \infty$, then Conditions Ri'$(\coneK^*)$, which 
is equivalent to Condition Ri, and Po'$(\coneK)$ hold by (iii) of Theorem~\ref{theorem:simple}.

Assertion (iv) also follows from assertion (iv) of Theorem~\ref{theorem:simple} 
since primal COP~\eqref{eq:COPprimal}, $\eta_d$, Po and Ri are identical or 
equivalent to  
P$(\coneK)$, $\zeta_d(\coneK^*)$, Po'$(\coneK)$ and Ri'$(\coneK^*)$, respectively. 
\qed

\smallskip

\begin{figure} 
\includegraphics[width=0.45\textwidth]{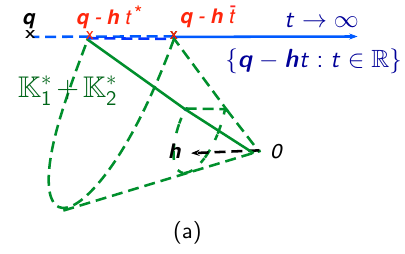}
\includegraphics[width=0.45\textwidth]{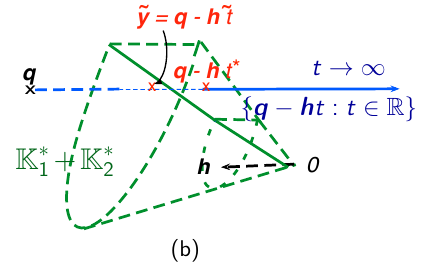}
\caption{Illustration for assertion (i) of Theorem~\ref{theorem:main0}} 
\end{figure}

Figure 2 illustrates assertion (i) of Theorem~\ref{theorem:main0}. 
In case (a),  neither Cl'$(\coneK_1^*+\coneK_2^*)$ nor 
Ri'($\coneK_1^*+\coneK_2^*$) is satisfied. 
The cone $\coneK_1^*+\coneK_2^*$ is almost open except the solid line that 
forms an extremal  
ray of the cone. The line $\{\q - \h t : t \in \Real\}$ touches the cone at the single point 
$\q-\h t^*$, but its line segment between $\q-\h t^*$ and $\q-\h \bar{t}$ is included 
in a face of 
its closure cl$(\coneK_1^*+\coneK_2^*) = (\coneK_1\cap\coneK)^*$. Hence 
$\zeta_d(\coneK_1^*+\coneK_2^*) =  t^* < \bar{t} =  \zeta_d((\coneK_1\cap\coneK_2)^*) 
= \zeta_p(\coneK_1\cap\coneK_2)$, which implies that a duality gap $\bar{t} - t^* > 0$ 
exists between 
primal COP P$(\coneK_1\cap\coneK_2)$ and its dual D$(\coneK_1^*+\coneK_2^*)$. 

The cone $\coneK_1^*+\coneK_2^*$ in case (b) is the same as in case (a), but 
 the line $\{\q - \h t : t \in \Real\}$ penetrates the interior of the cone, 
 so that Ri'($\coneK_1^*+\coneK_2^*$) is satisfied, where $\widetilde{\y} = \q - \h \widetilde{t} \in $ 
 int$(\coneK_1^*+\coneK_2^*)$. Both D$(\coneK_1^*+\coneK_2^*)$ and 
 D$((\coneK_1\cap\coneK_2)^*)$ share the common optimal value $t^*$ although 
 $\q - \h t^* \not\in \coneK_1^*+\coneK_2^*$ and D$(\coneK_1^*+\coneK_2^*)$ 
 has no feasible solution that attains the optimal value $t^*$. In this case, there is no duality 
 gap between primal COP P$(\coneK_1\cap\coneK_2)$ and its dual D$(\coneK_1^*+\coneK_2^*)$. 

\rema \label{Remark:main} 
Theorem~\ref{theorem:main0} could be proved by using \cite[Propositions 2.6 and 2.8]{SHAPIRO2001} and \cite[Theorems 5 and 7]{LUO1997}. But we need to convert 
COPs P$(\coneK_1\cap\coneK_2)$ and D$(\coneK_1^*+\coneK_2^*)$ into 
primal-dual pairs of COPs in different forms. See also Section 4. 
\erema

Necessary and/or sufficient conditions for the closedness of 
$\coneK_1^* + \coneK_2^*$ 
were throughly studied by Pataki \cite{PATAKI2007}. We refer to 
some of them to derive three important cases where Condition Cl holds. 
A closed convex cone $\coneJ \subset \spaceV$ is called nice 
if $\coneJ^* + F^{\perp}$ is closed for every face $F$ of~$\coneJ$. 
Polyhedral cones, positive semidefinite cones and doubly 
nonnegative cones are known to be nice. 
For every closed convex cone $\coneJ$ and $\x \in \coneJ$, let  
$\mbox{dir}(\x,\coneJ)$ denote the convex cone 
$\{\w : \x + t\w \in \coneJ  \ \mbox{for some } t > 0 \}$ and 
face$(\x,\coneJ)$ the smallest face of $\coneJ$ containing $\x$. 

\lemm \label{lemma:Pataki}
Let $\bar{\x} \in \mbox{relint}(\coneK_1\cap\coneK_2)$, 
$F_1 = \mbox{face}(\bar{\x},\coneK_1)$ and $F_2 = \mbox{face}(\bar{\x},\coneK_2)$. Then 
\begin{eqnarray}
\mbox{dir}(\bar{\x},\coneK_1)\cap\mbox{dir}(\bar{\x},\coneK_2)= 
(\mbox{cl dir}(\bar{\x},\coneK_1))\cap(\mbox{cl dir}(\bar{\x},\coneK_2))  
\label{eq:Pataki} 
\end{eqnarray}
is a necessary condition for $\coneK_1^* + \coneK_2^*$ to be closed. 
If in addition 
$ 
\coneK_1^* + F_1^{\perp} \ \mbox{ and } \coneK_2^* + F_2^{\perp} \ 
\mbox{are closed}  
$ 
--- in particular, if $\coneK_1$ and $\coneK_2$ are both nice --- 
then \eqref{eq:Pataki} is a sufficient condition. \vspace{-2mm} 
\elemm 
\proof{The assertion of the lemma is included in Theorem 5.1 of \cite{PATAKI2007}. In fact, \eqref{eq:Pataki} is given there as 
one of four equivalent 
necessary conditions for the closedness of $\coneK_1^* + \coneK_2^*$,
which are sufficient under the additional assumption mentioned above.\qed
}
 
\lemm \label{lemma:threeCases} 
Condition {\rm Cl}  holds ({\it i.e.}, $\coneK_1^* + \coneK_2^*$ is closed) if  one of the following conditions (i), (ii) and (iii) 
is satisfied.
\vspace{-2mm} 
\begin{description}
\item{(i)} $\coneK_1$ and $\coneK_2$ are polyhedral, {\it i.e.}, $\coneK_i = \{ \x \in \spaceV : 
\inprod{\a_i^j}{\x} \geq 0 \ (j=1,\ldots,\ell_i\}$ for some $\a_i^j \in \spaceV$ $(j=1,2,\ldots,\ell_i,i=1,2)$. 
\vspace{-2mm} 
\item{(ii)} 
There exists an $\widetilde{\x} \in ({\rm int}\coneK_1) \cap \coneK_2$.
\vspace{-2mm}
\item{(iii)} 
The cone $\coneK_1$ is nice, the cone $\coneK_2$ is polyhedral, and 
there exists an $\widetilde{\x} \in ({\rm relint}\coneK_1) \cap \coneK_2$.  
\vspace{-2mm} 
\end{description}
\elemm
\proof{
(i) In this case, we see that
\begin{eqnarray*}
(\coneK_1 \cap \coneK_2)^* & = & 
\left\{ \y= \sum_{j=1}^{\ell_1} \a^j_1u_j + \sum_{k=1}^{\ell_2} \a^k_2v_k  : 
\begin{array}{l}
u_j  \geq 0\ (j=1,\ldots,\ell_1), \\ 
v_k \geq 0 \ (k=1,\ldots,\ell_2) 
\end{array}
\right\} \ 
 = \ \coneK_1^* + \coneK_2^*. 
\end{eqnarray*}

(ii) Let $\widehat{\x} \in \mbox{relint}\coneK_2$. Since  $\widetilde{\x} \in ({\rm int}\coneK_1) \cap \coneK_2$, we can take a sufficiently small $\epsilon > 0$ such that
$
\bar{\x}  =  (1-\epsilon)\widetilde{\x} + \epsilon \widehat{\x} \in (\mbox{int}\coneK_1) \cap 
\mbox{relint}\coneK_2 \subset \mbox{relint}(\coneK_1\cap\coneK_2)
$. 
We see from $\bar{\x} \in \mbox{int}\coneK_1$ and $\bar{\x} \in \mbox{relint}\coneK_2$ 
that $F_i = $ face$(\bar{\x},\coneK_i) = \coneK_i$, which implies 
$\coneK_i^* + F_i^{\perp} = \coneK_i^*$,  and 
dir$(\bar{\x},\coneK_i) = \mbox{linspan}\coneK_i$ $(i=1,2)$. Since 
$\coneK_i^*$ and $\mbox{linspan}\coneK_i$ $(i=1,2)$ are closed, $\coneK_1^*+\coneK_2^*$ 
is closed by Lemma~\ref{lemma:Pataki}. 

(iii) Let $\hat{\x} \in \mbox{relint}(\coneK_1 \cap \coneK_2)$.  
Then 
$ 
\bar{\x}  =  ( \widetilde{\x} +  \hat{\x})/2 \in (\mbox{relint}\coneK_1) 
\cap \mbox{relint}(\coneK_1 \cap \coneK_2)
$ 
since $\widetilde{\x} \in \mbox{relint}\coneK_1$, 
$\hat{\x} \in \coneK_1$, 
$\hat{\x} \in \mbox{relint}(\coneK_1 \cap \coneK_2)$ 
 and $\widetilde{\x} \in \coneK_1 \cap \coneK_2$ hold. Since both cones 
 $\coneK_1$ and $\coneK_2$ are nice, it suffices to show~\eqref{eq:Pataki}.   
We see from 
 $\widetilde{\x} \in \mbox{relint}\coneK_1$ that dir$(\bar{\x},\coneK_1) = 
 $ linspan$\coneK_1$.
 Hence dir$(\bar{\x},\coneK_1) = $ \mbox{cl dir}$(\bar{\x},\coneK_1)$. Since $\coneK_2$ is 
 a polyhedral cone by the assumption, dir$(\bar{\x},\coneK_2)$ becomes polyhedral,  hence 
 dir$(\bar{\x},\coneK_2)  = $ cl dir$(\bar{\x},\coneK_2)$. 
 Therefore~\eqref{eq:Pataki} holds. 
\qed
}

The lemma below provides Slater type sufficient conditions for Conditions 
Ri$(\coneK_1^*+\coneK_2^*)$ and 
Po$(\coneK_1\cap\coneK_2)$. 

\lemm \label{lemma:sufficient} \mbox{ \ } \vspace{-2mm}
\begin{description}
\item{(i)} Conditions {\rm Ri} and {\rm Po} holds simultaneously if and 
only if 
there exists a $\widetilde{t} \in \Real$ such that  
$\q - \h \widetilde{t} \in \mbox{int}(\coneK_1^*+\coneK_2^*)$. 
\vspace{-2mm}
\item{(ii)} Assume that $\widetilde{\y}_1 \equiv \q - \h\widetilde{t} - \widetilde{\y}_2 
\in \mbox{int}(\coneK_1^*)$ for some $(\widetilde{t},\widetilde{\y}_2) 
\in \Real \times \coneK_2^*$. Then Conditions {\rm Ri} and {\rm Po} hold.
\vspace{-2mm} 
\end{description}
\elemm
\proof{
Assertion (i) follows from 
Lemmas~\ref{lemma:RiAndPo} and~\ref{lemma:MinkowskiSum}. 
We only prove assertion (ii). It suffices to show that 
$\widetilde{\y} \equiv \widetilde{\y}_1 + \widetilde{\y}_2 
= \q - \h \widetilde{t} \in \mbox{int}(\coneK_1^*+\coneK_2^*)$. 
Since $\widetilde{\y}_1  \in \mbox{int}\coneK_1^*$, 
we can take a positive number $\epsilon$ such that $\y_1 \in \coneK_1^*$ if 
$\|\y_1 - \widetilde{\y}_1\| \leq \epsilon$. 
Assume that $\y \in \spaceV$ and 
$\| \y - \widetilde{\y} \| \leq \epsilon$. Let $\y_1 =  \widetilde{\y}_1 + (\y-\widetilde{\y})$. Then  
$ \y = \y_1 + \widetilde{\y}_2$,  
$\|\y_1 - \widetilde{\y}_1\| = \|  \y-\widetilde{\y}\| \leq \epsilon$ 
(hence $\y_1 \in \coneK_1^*$), and $\widetilde{\y}_2 \in \coneK_2^*$ 
 by assumption. 
Therefore, $\y \in \coneK_1^* + \coneK_2^*$ and we have shown that $\widetilde{\y} \in \mbox{int}(\coneK_1^*+\coneK_2^*)$. 
\qed
}

%!TEX root = ./main.tex

\section{Strong duality in a symmetric primal-dual pair of COPs}

\label{section:symCOP}

Let $\spaceE_p$ and $\spaceE_d$ be finite dimensional vector spaces with inner products, 
which are both denoted by $\inprod{\cdot}{\cdot}$, and $\coneJ_p$ and $\coneJ_d$ closed 
convex cones of $\spaceE_p$ and $\spaceE_d$, respectively. Let $\b \in \spaceE_d$, 
$\cc\in \spaceE_p$, $\AC$ a linear map from 
$\spaceE_p$ into $\spaceE_d$. We denote the adjoint of $\AC$ by $\AC^*$. 
We consider the following primal-dual pair of COPs:
\begin{eqnarray}
\theta_p & = & \inf \left\{ \inprod{\cc}{\u} :  \u \in \coneJ_p, \ \AC \u -\b \in \coneJ_d \right\}.
\label{eq:COPsymP} \\ 
\theta_d & = & \sup \left\{ \inprod{\b}{\v} : \v \in \coneJ_d^* , \ \cc - \AC^*\v \in \coneJ_p^* \right\}. 
\label{eq:COPsymD} 
\end{eqnarray}
The following three convex cones play essential roles in characterizing strong 
duality of the primal-dual pair above in the subsequent discussion. 
\begin{eqnarray*}
\coneM_p 
& = & \left\{(\alpha ,\v-\AC\u) : \alpha \geq \inprod{\cc}{\u}, \ \u \in \coneJ_p, \ \v \in  \coneJ_d \right\}, \\
\coneN_p & = & \left\{\v-\AC\u : \u \in \coneJ_p, \ \v \in  \coneJ_d \right\}, \\ 
\coneN_p^* & = & 
\{\v \in \coneJ_d^*: -\AC^*\v \in \coneJ_p^* \}.  
\end{eqnarray*}

We derive the following result from Theorem~\ref{theorem:main0}. 
\theo \label{theorem:symmetricCOP}
Assume that $-\infty < \theta_p < \infty$ or $-\infty < \theta_d < \infty$. 
Then the assertions below hold. 
%{\bf \em assertions (i) and (iv) below}, and assertions (ii) and (iii) of Theorem~\ref{theorem:Shapiro} hold. 
\vspace{-2mm} 
\begin{description}
\item{(i)} If $\coneM_p$ is closed or $-\b \in {\rm relint}\coneN_p$, then 
$-\infty < \theta_p=\theta_d < \infty$ holds.\vspace{-2mm}
\item{(ii)} If $\coneM_p$ is closed, 
then primal COP~\eqref{eq:COPsymP} is solvable. \vspace{-2mm}
\item{(iii)} 
The set of optimal solutions of dual COP~\eqref{eq:COPsymD} is nonempty 
and bounded  if and only if $-\b \in {\rm int}\coneN_p$, 
which implies that $\coneN_p^*$ is pointed, holds. 
\vspace{-2mm}
\item{(iv)} 
The set of optimal  solutions of 
dual COP~\eqref{eq:COPsymD} is nonempty and unbounded if $\coneN_p^*$ is not pointed 
and  $-\b \in {\rm relint}\coneN_p$ holds.  
% then it is unbounded.
%}
\end{description}
\etheo

\rema 
$\coneM_p$ and $\coneN_p$ 
were originally introduced by Shapiro \cite{SHAPIRO2001} to establish 
strong duality between COPs~\eqref{eq:COPsymP} and ~\eqref{eq:COPsymD}. 
Theorem~\ref{theorem:symmetricCOP} is a slight modification (and extension) of 
\cite[Propositions 2.6 and 2.8]{SHAPIRO2001}. We note that  
Shapiro in \cite{SHAPIRO2001} dealt with the case where the spaces $\spaceE_p$ and 
$\spaceE_d$ can be a general vector 
(not necessarily finite dimensional) space. Theorem~\ref{theorem:symmetricCOP} 
could be proved by \cite[Propositions 2.6 and 2.8]{SHAPIRO2001} and 
\cite[Theorem 5]{LUO1997}. 
\erema

\noindent
As a corollary, we obtain:
\coro \label{coro:symmetricCOP}
If the set of optimal solutions (or the feasible region) of 
primal {\rm COP}~\eqref{eq:COPsymP} or 
dual 
{\rm COP}~\eqref{eq:COPsymD} is nonempty and bounded, then $-\infty < \theta_p=\theta_d < \infty$. \vspace{-3mm}
\ecoro
\proof{
Assume that the set of optimal solutions  of 
dual COP~\eqref{eq:COPsymD} is nonempty and bounded. Then, 
by the assumption and (iii) of Theorem~\ref{theorem:symmetricCOP}, we see that 
$-\infty < \theta_d < \infty$ and $-\b \in {\rm int}\coneN_p$ hold. By (i) of Theorem~\ref{theorem:symmetricCOP}, we obtain that $-\infty < \theta_p=\theta_d < \infty$. 
Since the primal-dual pair of COPs~\eqref{eq:COPsymP} and~\eqref{eq:COPsymD} is symmetric, 
we can derive the same conclusion even if we regard COP~\eqref{eq:COPsymD} as 
primal and COP~\eqref{eq:COPsymP} as dual. See the  discussions in 
Sections~\ref{subsection:symCOPpp} and~\ref{subsection:symCOPpd}. This corollary 
could also be obtained from Corollary~\ref{coro:significant}
\qed
}

\smallskip

We convert primal COP~\eqref{eq:COPsymP} to COP~\eqref{eq:COPprimal} in 
Sections~\ref{subsection:symCOPpp},  and dual COP~\eqref{eq:COPsymD} to COP~\eqref{eq:COPprimal}
in Section~\ref{subsection:symCOPpd}. 
The proof of Theorem~\ref{theorem:symmetricCOP} is stated 
in Section~\ref{subsection:proof} where the latter conversion given in 
Section~\ref{subsection:symCOPpd} is utilized. 
Although Sections~\ref{subsection:symCOPpp} is not relevant 
to the proof of Theorem~\ref{theorem:symmetricCOP}, it  may clarify the discussion in Section~\ref{subsection:symCOPpd} and the dual counter part of Theorem~\ref{theorem:symmetricCOP}. 
The discussion in Section ~\ref{subsection:symCOPpp} will help the reader to derive such  counter part. 
Table~2 summarizes the duality 
results shown in Theorems~\ref{theorem:symmetricCOP}, Corollary~\ref{coro:symmetricCOP} and the  results on their duals % counter parts 
(see Section~\ref{subsection:symCOPpp}). 

%\afterpage{
\begin{table}[h!]
\begin{center}
%\scriptsize{
\begin{tabular}{|c||c|c|}
\hline
%                                            &
\multicolumn{1}{|c||}{Sufficient conditions}  & \multicolumn{2}{|c|}
{The sets of primal and dual optimal solutions}  \\       
\multicolumn{1}{|c||}{for strong duality $\theta_p=\theta_d$} & pSol & dSol \\ 
\hline
%\hline
$\coneM_p$ : closed & nonempty & ? \\
$\coneN_p^*$ : pointed \& $-\b\in{\rm relint}\coneN_p$ & ? & nonempty \& bounded \\ 
$\coneN_p^*$ : not pointed \& $-\b\in{\rm relint}\coneN_p$ & ? & nonempty \& unbounded  \\ 
dSol : nonempty \& bounded & ? & nonempty \& bounded \\
\hline
$\coneM_d$ : closed & ? & nonempty \\
$\coneN_d^*$ : pointed \& $ \cc\in{\rm relint}\coneN_d$ & nonempty \& bounded 
& ? \\ 
$\coneN_d^*$ : not pointed \& $\cc\in{\rm relint}\coneN_d$ & nonempty \& unbounded & ? \\ 
pSol : nonempty \& bounded & nonempty \& bounded  & ? \\
\hline                                             
\end{tabular}	
\end{center}
\caption{
%\textcolor{red}{
A summary of strong duality $-\infty < \theta_p=\theta_d < \infty$, 
solvability of primal COP~\eqref{eq:COPsymP} 
and solvability of dual COP~\eqref{eq:COPsymD}. 
%Throughout the table,
For all cases, we assume that $-\infty < \theta_p < \infty$ or 
$-\infty < \theta_d < \infty$ holds. 
``pSol" and ``dSol" denote 
the sets of optimal solutions of COP~\eqref{eq:COPsymP} and COP~\eqref{eq:COPsymD}, respectively. If $\coneN_p^*$ ($\coneN_d^*$) is pointed, 
then ${\rm relint}\coneN_p = {\rm int}\coneN_p$ (${\rm relint}\coneN_d = {\rm int}\coneN_d$).
$\coneM_p 
= \left\{(\alpha ,\v-\AC\u) : \alpha \geq \inprod{\cc}{\u}, \ \u \in \coneJ_p, \ \v \in  \coneJ_d \right\}$, 
$\coneN_p = \left\{\v-\AC\u : \u \in \coneJ_p, \ \v \in  \coneJ_d \right\}$, 
$\coneN_p^* = 
\{\v \in \coneJ_d^*: -\AC^*\v \in \coneJ_p^* \}$, 
$\coneM_d = \left\{ (\beta,\u + \AC^*\v) : \beta \leq \inprod{\b}{\v}, \ \u \in \coneJ_p^*, \v \in \coneJ_d^* \right\}$, 
$\coneN_d = \left\{ \u + \AC^*\v :  \u \in \coneJ_p^*, \v \in \coneJ_d^* \right\}$, $\coneN_d^* = \{\u \in \coneJ_p: \AC\u \in \coneJ_d\}$. 
``?'' means that the existence of an optimal solution is not guaranteed.
%}
}
\end{table}
%}

\subsection{Conversion of the primal-dual pair of COP~\eqref{eq:COPsymP} and COP~\eqref{eq:COPsymD} to the 
primal-dual pair of COP~\eqref{eq:COPprimal} and COP~\eqref{eq:COPdual}}

\label{subsection:symCOPpp}

Let 
\begin{eqnarray}
\left.
\begin{array}{lcl}
\spaceV & = & \Real \times \spaceE_p = \left\{ (x_0,\u) : x_0 \in \Real,\ \u \in \spaceE_p\right\}, \\
\coneK_1 & = & \Real_+ \times \coneJ_p, \ 
\coneK_2 = \left\{ (x_0,\u) \in \spaceV : x_0 \in \Real,\ \AC \u - \b x_0 \in \coneJ_d \right\}, \\  
\q & = & (0,\cc) \in \spaceV, \ \h = (1,\0) \in \coneK_1^*. 
\end{array}
\right\} \label{eq:DefSymP-Primal}
\end{eqnarray}
It should be noted that the symbols $\spaceV$, $\coneK_1$, $\coneK_2$, $\q$ and $\h$ above are defined locally 
in this section. The same symbols will be also defined locally in Section~\ref{subsection:symCOPpd} with different meanings. 
As the inner product of $\x = (x_0,\u) \in \spaceV$, $\y = (t_0,\w) \in \spaceV$, we 
define $\inprod{\x}{\y} = x_0t_0 + \inprod{\u}{\w}$.  
Then, 
\begin{eqnarray*}
& & \coneK_1^* = \Real_+ \times \coneJ_p^*, \ 
\coneK_2^* = \left\{(\inprod{-\b}{\v},\AC^*\v) : \v \in \coneJ_d^* \right\}, \\
& &  
\u \in \coneJ_p, \ \AC \u -\b \in \coneJ_d \ 
\mbox{if and only if } 
\x  = (x_0,\u) \in \coneK_1, \ \x \in \coneK_2, \ \inprod{\h}{\x} = 1.
\end{eqnarray*}
Thus, we can rewrite COP~\eqref{eq:COPsymP} as COPs~\eqref{eq:COPprimal}. 
We also see that 
\begin{eqnarray*}
%%%%%%%%%%
& & 
 \q - \h t - \y_2 = (0,\cc)  - (t,0) - (\inprod{-\b}{\v},\AC^*\v) \\
&& \mbox
{ \ } \hspace{20.5mm}  = 
(\inprod{\b}{\v} - t,\cc -  \AC^*\v) \\ 
& & { \ } \hspace{28mm} \mbox{if } \y_2 = (\inprod{-\b}{\v},\AC^*\v) \in \coneK_2^* \
\mbox{ for some } \v \in \coneJ_d^*; \ \mbox{hence}\\
%%%%%%%%%%
& & 
\inprod{\b}{\v} \geq t  \ 
\mbox{and } \cc - \AC^*\v \in \coneJ_p^* \ \mbox{for some } \v \in \coneJ_d^* \ 
\mbox{if and only if } 
\q - \h t \in \coneK_1^* + \coneK_2^*. 
\end{eqnarray*}
Therefore,  COP~\eqref{eq:COPsymD} can be rewritten as COP~\eqref{eq:COPdual}. 
In this case, $\coneK_1^* + \coneK_2^*$ is represented as 
$\left\{ (-\inprod{\b}{\v}+s,\AC^*\v+\u) : s \in \Real_+, \ \u \in \coneJ_p^*, \v \in \coneJ_d^* \right\}$. 
Obviously, $\coneK_1^* + \coneK_2^*$ is closed if and only if 
\begin{eqnarray*}
\coneM_d & = & \left\{ (\beta,\u + \AC^*\v) : \beta \leq \inprod{\b}{\v}, \ \u \in \coneJ_p^*, \v \in \coneJ_d^* \right\}
\end{eqnarray*} 
is closed. 
The lemma below 
converts Conditions Ri ($= $ {\rm Ri'}$(\coneK_1^*+\coneK_2^*)$) 
and Po ($= $Po'$(\coneK_1\cap\coneK_2)$) 
on COP~\eqref{eq:COPdual} 
to the corresponding conditions on COP~\eqref{eq:COPsymD}, respecctively.

\lemm \label{lemma:TdForCOPsymD}
Assume that $\spaceV$, $\coneK_1$, $\coneK_2$, $\q$ and $\h$ are given in~\eqref{eq:DefSymP-Primal}.\vspace{-3mm}
\begin{description}
\item{(i)} Condition {\rm Ri'}$(\coneK_1^*+\coneK_2^*)$ 
holds if and only if $\cc \in  {\rm relint}\coneN_d$, where \vspace{-2mm}
\begin{eqnarray*}
\coneN_d = \left\{ \u + \AC^*\v :  \u \in \coneJ_p^*, \v \in \coneJ_d^* \right\}.  
\end{eqnarray*}
\vspace{-10mm}
\item{(ii)} $\coneK_1\cap\coneK_2$ is pointed if and only if 
$\coneN_d^*$ is pointed.
\end{description}
\vspace{-6mm}
\elemm 
\proof{
``only if part'' of (i):  Assume that 
$\q - \h \widetilde{t} = (-\widetilde{t},\cc)$ lies in relint$(\coneK_1^*+\coneK_2^*)$. 
Then there exists 
$(\widetilde{s},\widetilde{\u},\widetilde{\v}) \in \Real_+\times \coneJ_p^* \times \coneJ_d^*$ such that 
\begin{eqnarray*}
(-\widetilde{t},\cc) & = & (\inprod{-\b}{\widetilde{\v}} + \widetilde{s}, \widetilde{\u} + \AC^* \widetilde{\v}) \\ 
                   & \in & \mbox{relint}  (\coneK_1^*+\coneK_2^*) \\ 
                   & = & \mbox{relint} \left\{(\inprod{-\b}{\v} + s, \u + \AC^* \v) : s \in \Real_+, \ \u \in \coneJ_p^*, \ 
                   \v \in \coneJ_d^* \right\},   
\end{eqnarray*} 
which implies $\cc = \widetilde{\u} + \AC^* \widetilde{\v} \in $ relint$\coneN_d$. 

``if part'' of (i): Assume that $\cc = \widetilde{\u} + \AC^* \widetilde{\v} \in $ relint$\coneN_d$. 
If we take $\widetilde{s} >0$ and let $\widetilde{t} = - (\inprod{-\b}{\widetilde{\v}} + \widetilde{s})$, 
then $\q - \h \widetilde{t} = (-\widetilde{t},\cc) = (\inprod{-\b}{\widetilde{\v}} + \widetilde{s}, \widetilde{\u} + \AC^* \widetilde{\v})\in$ relint$(\coneK_1^*+\coneK_2^*)$. 

(ii): We prove contraposition of the assertion. To prove ``if part', assume that 
$\coneK_1\cap\coneK_2$ is not pointed. Then there exists a nonzero $(x_0,\u) \in \spaceV$ such that 
\begin{eqnarray}
(x_0,\u) \in \coneK_1\cap\coneK_2 \ \mbox{ and } -(x_0,\u) \in \coneK_1\cap\coneK_2. \label{eq:notPointedCOPprimal}
\end{eqnarray}
It follows that $x_0 = 0$, $\u \not= \0$, 
\begin{eqnarray}
\u \in \{\u \in \coneJ_p: \AC\u \in \coneJ_d\} \ \mbox{ and }
-\u \in \{\u \in \coneJ_p: \AC\u \in \coneJ_d\}. \label{eq:notPointedForCOPsymD}
\end{eqnarray}
This implies $\{\u \in \coneJ_p: \AC\u \in \coneJ_d\}$ is not pointed. 
Now, assume that $\{\u \in \coneJ_p: \AC\u \in \coneJ_d\}$ is not pointed to prove ``only if part''. Then there exists a nonzero $\u$ 
satisfying~\eqref{eq:notPointedForCOPsymD}. Let $x_0=0$. Then 
$\0 \not= (x_0,\u) \in \spaceV$ satisfies~\eqref{eq:notPointedCOPprimal}, 
which implies that $\coneK_1\cap\coneK_2$ is pointed. We can easily verify 
the identity $\{\u \in \coneJ_p: \AC\u \in \coneJ_d\} = \coneN_d^*$ 
as the the identity $\{\v \in \coneJ_d^* : -\AC^*\v \in \coneJ_p^* \} = 
\coneN_p^*$. 
%  whose proof is given in Section~\ref{section:Np}.
\qed 
}

\subsection{Conversion of the primal-dual pair of COP~\eqref{eq:COPsymP} and COP~\eqref{eq:COPsymD} to the 
dual-primal pair of COP~\eqref{eq:COPdual} and COP~\eqref{eq:COPprimal}}

\label{subsection:symCOPpd}

Since the primal-dual pair of COPs~\eqref{eq:COPsymP} and~\eqref{eq:COPsymD} is symmetric, 
we can  interchange their role as follows:
\begin{eqnarray}
\eta_p & = & -\theta_d  = 
\inf \left\{ \inprod{-\b}{\v} : \v \in \coneJ_d^* , \ \cc - \AC^*\v \in \coneJ_p^* \right\}\nonumber  \\ 
& =  &\inf \left\{ \inprod{\widetilde{\b}}{\v} : \v \in \coneJ_d^* , \ -\widetilde{\cc} - \AC^*\v \in \coneJ_p^* \right\}. 
\label{eq:COPsymP-Dual} \\ 
\eta_d & = & -\theta_p =  \sup \left\{ \inprod{-\cc}{\u} :  \u \in \coneJ_p, \ \AC \u -\b \in \coneJ_d \right\}\nonumber \\ 
& = & \sup \left\{ \inprod{\widetilde{\cc}}{\u} :  \u \in \coneJ_p, \ \AC \u +\widetilde{\b} \in \coneJ_d \right\}.
\label{eq:COPsymD-Primal}
\end{eqnarray}
Here $\widetilde{\b} = -\b$ and $\widetilde{\cc} = -\cc$. Now we regard COP~\eqref{eq:COPsymP-Dual} 
induced  from dual COP~\eqref{eq:COPsymD} as primal,  and 
COP~\eqref{eq:COPsymD-Primal} induced  from primal COP~\eqref{eq:COPsymP} as dual. Let 
\begin{eqnarray}
\left.
\begin{array}{lcl} 
\spaceV & = & \Real \times \spaceE_d = \left\{ (x_0,\v) : x_0 \in \Real,\ \v \in \spaceE_d\right\},\\[2pt]
\coneK_1 & = & \Real_+ \times \coneJ_d^*, \ 
\coneK_2 = \left\{ (x_0,\v) \in \spaceV :  - \widetilde{\cc} x_0 -\AC^* \v  \in \coneJ_p^* \right\}, \\[2pt]  
\q & = & (0,\widetilde{\b}) \in \spaceV, \ \h = (1,\0) \in \coneK_1^*. 
\end{array}
\right\}
\label{eq:DefSymP-Dual} 
\end{eqnarray}
Then we can similarly show as in Section~\ref{subsection:symCOPpp} 
that COPs~\eqref{eq:COPsymP-Dual} and~\eqref{eq:COPsymD-Primal} are equivalent to COPs~\eqref{eq:COPprimal} and~\eqref{eq:COPdual}, respectively. 
As a result, the pair of dual COP~\eqref{eq:COPsymD}
and primal COPs~\eqref{eq:COPsymP} is equivalently reformulated as the primal-dual pair of COPs~\eqref{eq:COPprimal} and COPs~\eqref{eq:COPdual}
with $\spaceV$, $\coneK_1$, $\coneK_2$, $\q$ and $\h$ given in~\eqref{eq:DefSymP-Dual}.
We also see that  
\begin{eqnarray*}
\coneK_1^* + \coneK_2^*
& = & \left\{(\alpha ,\v-\AC\u) : \alpha \geq \inprod{\cc}{\u}, \ \u \in \coneJ_p, \ \v \in  \coneJ_d \right\}
  = \coneM_p. 
\end{eqnarray*}
The following lemma can be  proved similarly to Lemma~\ref{lemma:TdForCOPsymD}. 
\lemm \label{lemma:TdForCOPsymP} 
Assume that $\spaceV$, $\coneK_1$, $\coneK_2$, $\q$ and $\h$ are given in~\eqref{eq:DefSymP-Dual}. \vspace{-2mm}
\begin{description}
\item{(i)} Condition {\rm Ri'}$(\coneK_1^*+\coneK_2^*)$ 
holds if and only if $-\b \in {\rm relint}\coneN_p$.
\vspace{-2mm}
\item{(ii)} $\coneK_1\cap\coneK_2$ is pointed if and only if 
$\coneN_p^*$ is pointed.
\end{description}
\elemm 

\subsection{Proof of Theorem~\ref{theorem:symmetricCOP}} 

\label{subsection:proof}

We reformulate the pair of dual COP~\eqref{eq:COPsymD}
and primal COP~\eqref{eq:COPsymP}  
as the pair of primal COP~\eqref{eq:COPprimal} 
and dual COP~\eqref{eq:COPdual} with $\spaceV$, $\coneK_1$, $\coneK_2$, $\q$ and $\h$ given in~\eqref{eq:DefSymP-Dual} as discussed in the previous section, and 
apply Theorem~\ref{theorem:main0} to the reformulated pair. We know that 
$\eta_p = -\theta_d$, $\eta_d = -\theta_p$, 
%\red{
Condition Cl'$(\coneK_1^*+\coneK_2^*)$ is equivalent to the closedness of 
$\coneM_p$, Condition Ri'$(\coneK_1^*+\coneK_2^*)$ is equivalent to 
$-\b \in $ relint$\coneN_p$, and $\coneK_1\cap\coneK_2$ is pointed if and 
only if $\{\v \in \coneJ_d^*: -\AC^*\v \in \coneJ_p^* \} = \coneN_p^*$ is 
pointed. 
We also know by assertion (i) 
of Lemma~\ref{lemma:sufficient} that $-\b \in $ int$\coneN_p$ if and only if 
Conditions Ri and Po are satisfied. 
Therefore, all assertions of Theorem~\ref{theorem:symmetricCOP} follow from 
Theorem~\ref{theorem:main0}. 
\qed

%
%!TEX root = ./main.tex

\section{Concluding remarks}

\label{section:concludingRemarks}
 
 In  \cite{KIM2013}, Kim, Kojima and Toh formulated 
a DNN relaxation of a class of linearly 
 constrained QOPs in nonnegative and binary variables as a COP of the form
\begin{eqnarray}
\varphi_p & = & \inf \left\{ \inprod{\q}{\x} : \x \in \coneK_1, \ \inprod{\h}{\x} = 1, \ \inprod{\h_1}{\x} = 0
\right\}, \label{eq:QopDNNprimal} 
\end{eqnarray} 
where $\coneK_1$ is a nice cone (the intersection of the doubly nonnegative 
cone and a linear subspace) in the space of symmetric matrices, 
$\0 \not=\h \in \coneK_1^*$ and $\h^{1} \in \coneK_1^*$. 
The dual of COP~\eqref{eq:QopDNNprimal} can be written as 
\begin{eqnarray}
\varphi_d & = & \sup \left\{ t : \q - \h t - \h s \in \coneK_1^* \right\}. 
\label{eq:QopDNNdual} 
\end{eqnarray}
Under the assumption that ensures the feasible region of COP~\eqref{eq:QopDNNprimal} is 
nonempty and bounded, they proved the strong duality that $-\infty < \varphi_p =  \varphi_d < \infty$ 
\cite[Lemma 3]{KIM2013} (see also \cite[Theorem 2.6]{ARIMA2018} for 
the same assertion under a weaker 
assumption). By applying a Lagrangian relaxation to COPs~\eqref{eq:QopDNNprimal} 
and~\eqref{eq:QopDNNdual},  they induced the primal-dual pair of COPs of the form~\eqref{eq:COPprimal}
and the form~\eqref{eq:COPdual}, which has no duality gap, 
with a Lagrangian multiplier parameter $\lambda$ associated with 
the constraint $\inprod{\h_1}{\x} = 1$. The strong duality equality 
$\varphi_p = \varphi_d$ was derived by taking the limit of their optimal values with no gap as 
$\lambda \rightarrow \infty$. This result was used in the development of the 
Newton-bracketing method whose convergence is quadratic 
 for solving the 
Lagrangian-DNN relaxation of the aforementioned class of QOPs in \cite{KIM2021}. 
See \cite{ARIMA2018,KIM2021} for more details. 

The primal-dual pair of COPs~\eqref{eq:QopDNNprimal} and~\eqref{eq:QopDNNdual} can be reformulated 
as the primal-dual pair of COPs~\eqref{eq:COPprimal} and~\eqref{eq:COPdual} 
with $\coneK_2 = \{ \x \in \spaceV : 
\inprod{\h_1}{\x} = 0 \}$. Under their assumption, 
the set of optimal solutions of primal COP~\eqref{eq:COPprimal}
is nonempty and bounded. Hence, their strong duality 
can be derived directly from Corollary~\ref{coro:significant} (or Theorem~\ref{theorem:main0})    
without relying on the Lagrangian relaxation of COP~\eqref{eq:QopDNNprimal}.

\bibliographystyle{plain}
\bibliography{./enhFOM}

%%%%%

\end{document}